\documentclass[12pt]{article}
\usepackage{amsmath}
\usepackage{amsfonts}
\usepackage{amssymb}
\usepackage{rotating}
\usepackage{mathdots}
\usepackage[usenames,dvipsnames]{pstricks}
\usepackage{epsfig}
\usepackage{pst-grad}
\usepackage{pst-plot}
\usepackage{slashbox}
\usepackage{ctable}

\def\tref#1{Theorem~\ref{#1}}
\def\eref#1{(\ref{#1})}
\def\sref#1{\S\ref{#1}}
\def\va{{\rm\bf a}}
\def\vb{{\rm\bf b}}
\def\vc{{\rm\bf c}}
\def\vd{{\rm\bf d}}
\def\ve{{\rm\bf e}}
\def\vx{{\rm\bf x}}
\def\vy{{\rm\bf y}}
\def\vz{{\rm\bf z}}
\def\vv{{\rm\bf v}}
\def\vw{{\rm\bf w}}
\def\vr{{\rm\bf r}}
\def\vs{{\rm\bf s}}
\def\vu{{\rm\bf u}}
\def\vt{{\rm\bf t}}
\def\vm{{\rm\bf m}}
\def\vo{{\rm\bf o}}

\newtheorem{theorem}{Theorem}[section]
\newtheorem{proposition}[theorem]{Proposition}

\begin{document}

\author{Gi-Sang Cheon$^a$\footnote{This work was supported by the National
Research Foundation of Korea Grant funded by the Korean Government
(NRF-2012-007648)}, Hana Kim$^b$\footnote{This research was
supported by Basic Science Research Program through the National
Research Foundation of Korea(NRF) funded by the Ministry of
Education, Science and Technology
(2013R1A6A3A03024342)} and Louis W. Shapiro$^c$\\
{\footnotesize $^a$ \textit{Department of Mathematics, Sungkyunkwan
University, Suwon 440-746, Rep. of Korea}}\\
{\footnotesize gscheon@skku.edu}\\
{\footnotesize $^b$ \textit{National Institute for Mathematical
Sciences, 70 Yuseong-daero, 1689 beon-gil,
Yuseong-gu,}}\\{\footnotesize\textit{Daejeon 305-811, Rep. of
Korea}}\\
{\footnotesize hkim@nims.re.kr}\\
{\footnotesize $^c$ \textit{Department of Mathematics, Howard
University,
Washington, DC 20059, USA}}\\
{\footnotesize lshapiro@howard.edu}}
\title{Mutation effects in ordered trees}
\date{}\maketitle

\begin{abstract}
A mutation will affect an individual and some or all of its
descendants. In this paper, we investigate ordered trees with a
distinguished vertex called the mutator. We describe various
mutations in ordered trees, and find the generating functions for
statistics concerning trees with those mutations. The examples give
new interpretations to several known sequences and also introduce
many new sequences and their combinatorial interpretations.
\end{abstract}

\vskip1pc \noindent\textit{AMS classifications}: Primary: 05A15;
secondary: 05C05

\noindent\textit{Key words}: mutation, ordered trees, short lived
mutation, toggle tree, right path tree.
\bigskip

\section{Introduction}

A mutation is a change in the genome of an organism as well as a
genotype that exhibits high rates of mutation. It can result in
several different types of change in the nucleotide sequences. For
instance, mutations in genes can either have no effect, alter the
product of a gene, or prevent the gene from functioning properly or
completely. These phenomena may be reflected in a family tree which
is an example of ordered tree where its subtrees are usually ordered
by date of birth but could be ordered by some other attribute such
as height.

In the present paper, we consider ordered trees with one
distinguished vertex called the \textit{mutator}. The vertices
changed by the mutator are said to be of a \textit{new type}. We
include the mutator itself as being of the new type. Ordered trees
with a mutation may reflect many biological or social structures
where changes occur. There are an abundant literature on
applications. For instance, if the mutator represents a genetic
mutation in a family then the new type vertices are those carrying
this mutation. If an ordered tree represents a river network then a
mutator could be a spot where pollution has been detected and the
vertices above it could be the possible source of the pollution.

It is well-known that the number of ordered trees with $n$ edges is
counted by the $n$th Catalan number $C_n={1\over n+1}{2n\choose n}$,
and its generating function is
$C=\sum_{n\ge0}C_nz^n={1-\sqrt{1-4z}\over2z}$ resulting from
$C=1+zC^2={1\over1-zC}$. In particular, the number of the ordered
trees with a distinguished vertex is counted by the $n$th central
binomial coefficient $B_n={2n\choose n}$ with
$B=\sum_{n\ge0}B_nz^n={1\over\sqrt{1-4z}}$. A key fact directly
obtained from the uplift principle (see Proposition \ref{uplift})
provides the generating function $B/C$ for the ordered trees with a
distinguished leaf (i.e. terminal vertex). \vskip.5pc

The purpose of this paper is to investigate five kinds of ordered
trees with mutations according to some conditions on the children of
the mutator. In particular, we postulate that such conditions are
given or can be explained in terms of generating functions. We then
enumerate the numbers of such trees and their vertices as well as
vertices of each new type. The asymptotic behavior of those numbers
will be also discussed. Finally in Section 3 we change from ordered
trees to complete binary trees.

\section{Ordered trees with various mutations}

Throughout this section, a tree means an ordered tree. Let us begin
with trees affected by a mutation. By the children we shall mean the
vertices directly connected to the mutator, and by the descendants
we shall mean all subsequent vertices above the mutator. There are a
variety of types of conditions that we could set for the children or
the descendants of the mutator.

At one extreme the mutation causes sterility. In this case the
generating function for the number of trees is $T_M:=L=B/C$ and we
have only one new type vertex, the mutator per tree. At the other
extreme all the descendants of the mutator are of the new type. As
an intermediate case once a new type child appears all subsequent
descendants are new type. This could model a case where the mutator
is a person who moves to a new country. Other possibilities include
exactly one new and sterile child. We also consider what happens if
the new type vertices are all on the right most branch or on the
right most path.

There are many further variations possible. In one direction we
could look at different kinds of trees such as complete (or
incomplete) binary trees, Motzkin 0-1-2 trees, even trees, Riordan
trees, and spoiled child trees. A second direction would allow a
variety of new types, and a third would allow more than one mutator.
All of mutation possibilities can be considered with these
variations, but to keep this paper focused we discuss these only
briefly.

In this section, we consider five kinds of trees with mutations
arising from different conditions on the children or descendants. We
compute the number of vertices of each new type by the following
procedure. First we find the generating function for such trees,
designate this generating function as $T_M$. Since each tree with
$n$ edges has $n+1$ vertices, the generating function for the
vertices is the derivative $(zT_M)^\prime$. We then find the
generating function $V_N$ for vertices of the new type. In this
step, the uplift principle can be used for transferring results
established at the root to an arbitrary vertex. We note that a
mutator does not change the conditions on its children wherever it
appears as is usual.

\begin{proposition}\label{uplift} {\rm (The uplift principle, \cite{BCheon1})}
First, find the generating function for whatever is being counted at
the root. Then uplift the result at the root to an arbitrary vertex
by multiplying by the leaf generating function $L=B/C$.
\end{proposition}

Next, we compute the proportion of new type vertices among all the
vertices. For this step we use a few asymptotic results depending on
Stirling's approximation or the ratio test. These are
\begin{itemize}
\item[(i)] ${2n\choose n}\sim{4^n\over\sqrt{\pi n}}$ or
$4^n\sim{2n\choose n}\sqrt{\pi n}$, ${2n\choose n}\sim4{2n-2\choose
n-1}$;
\item[(ii)] $C_n\sim4C_{n-1}$, $B_n\sim4B_{n-1}$;

\end{itemize}
We also adopt a singularity analysis for few complicated cases. In
addition here are some other facts we will use:
\begin{itemize}
\item[(i)] $L=B/C$, $C^\prime=BC^2,\;\; B^\prime=2B^3$;
\item[(ii)] $B=1+2zBC={1\over1-2zC}={C\over1-zC^2}\;\;\hbox{so that}\;\;{B-1\over2}=zBC$;
\item[(iii)] $[z^n]C^s={s\over
2n+s}{2n+s\choose n}\;\;\hbox{and}\;\;[z^n]BC^s={2n+s\choose n}$.
\end{itemize}
where $[z^n]$ is the coefficient extraction operator.\vskip.5pc

Let us now describe first example of ordered trees with a mutation
under the variation of an extreme case. Here the rightmost child of
the mutator is of the new type. The mutator has no more children of
the new type and the child of the new type has no children. In this
sense, we call such mutation the \textit{short lived mutation}. As
an example if a male donkey mates with a female horse, it stops
reproducing and the child, a mule, is sterile. The horse,
traumatized by the experience, has no more offspring.

Every such tree has exactly two vertices of the new type so the
question of interest is the number of such trees. The generating
function is
\begin{eqnarray*}
L\cdot C\cdot z={B\over C}\cdot C\cdot
z=zB={z\over\sqrt{1-4z}}=z+2z^2+6z^3+20z^4+\cdots\quad (A000984).
\end{eqnarray*}
To illustrate with $n=3$ edges, we consider the ordered trees with a
short lived mutator and three edges. There are $6$ trees as shown in
Figure 1. From now on, `x' denotes the root and the mutator is
circled. Along with the mutator, we will mark the edges above the
vertices of the new type as this is easier to see.
\begin{figure}[h]
\begin{center}
\scalebox{0.6} {
\begin{pspicture}(0,-1.98)(12.032532,2.0)
\psline[linewidth=0.04cm](0.26,1.84)(0.26,-1.76)
\psline[linewidth=0.04cm](6.0824585,-1.781012)(6.657208,-0.4587174)
\psline[linewidth=0.04cm](6.0824585,-1.781012)(5.477217,-0.463407)
\psdots[dotsize=0.24](0.26,1.86) \psdots[dotsize=0.24](0.26,0.66)
\psline[linewidth=0.08cm](0.12,-1.58)(0.4,-1.88)
\psline[linewidth=0.08cm](0.4,-1.58)(0.12,-1.88)
\psline[linewidth=0.04cm](1.86,-0.52)(1.86,-1.72)
\psdots[dotsize=0.24](1.86,-0.52)
\psline[linewidth=0.08cm](1.7,-1.54)(1.98,-1.84)
\psline[linewidth=0.08cm](1.98,-1.54)(1.7,-1.84)
\psdots[dotsize=0.24](5.46,-0.52) \psdots[dotsize=0.24](6.64,-0.52)
\psline[linewidth=0.08cm](5.92,-1.56)(6.2,-1.86)
\psline[linewidth=0.08cm](6.2,-1.56)(5.92,-1.86)
\psline[linewidth=0.08cm](3.7,-1.54)(3.98,-1.84)
\psline[linewidth=0.08cm](3.98,-1.54)(3.7,-1.84)
\psline[linewidth=0.04cm](10.88,-1.721012)(10.88,-0.52)
\psdots[dotsize=0.24](10.88,-0.56)
\pscircle[linewidth=0.03,dimen=outer](0.26,0.66){0.26}
\pscircle[linewidth=0.03,dimen=outer](1.84,-0.52){0.26}
\psdots[dotsize=0.24](0.26,-0.54)
\psline[linewidth=0.04cm](0.14,1.7)(0.38,1.7)
\psline[linewidth=0.04cm](0.14,1.56)(0.38,1.56)
\psline[linewidth=0.04cm](0.14,1.4)(0.38,1.4)
\psline[linewidth=0.04cm](0.14,1.26)(0.38,1.26)
\psline[linewidth=0.04cm](0.14,1.14)(0.38,1.14)
\psline[linewidth=0.04cm](0.14,1.0)(0.38,1.0)
\psline[linewidth=0.04cm](0.14,0.86)(0.38,0.86)
\pscircle[linewidth=0.03,dimen=outer](5.44,-0.52){0.26}
\psline[linewidth=0.04cm](5.32,0.52)(5.56,0.52)
\psline[linewidth=0.04cm](5.32,0.38)(5.56,0.38)
\psline[linewidth=0.04cm](5.32,0.22)(5.56,0.22)
\psline[linewidth=0.04cm](5.32,0.08)(5.56,0.08)
\psline[linewidth=0.04cm](5.32,-0.04)(5.56,-0.04)
\psline[linewidth=0.04cm](5.32,-0.18)(5.56,-0.18)
\psline[linewidth=0.04cm](5.32,-0.32)(5.56,-0.32)
\psline[linewidth=0.04cm](5.44,0.64)(5.44,-0.56)
\psdots[dotsize=0.24](5.46,0.68)
\psline[linewidth=0.04cm](1.8824584,-0.581012)(2.457208,0.7412826)
\psline[linewidth=0.04cm](1.8824584,-0.581012)(1.2772171,0.736593)
\psdots[dotsize=0.24](1.26,0.68) \psdots[dotsize=0.24](2.44,0.68)
\psline[linewidth=0.04cm](2.2198203,0.53498775)(2.4436145,0.44829217)
\psline[linewidth=0.04cm](2.1692476,0.40444112)(2.393042,0.3177455)
\psline[linewidth=0.04cm](2.1114507,0.25524494)(2.335245,0.16854934)
\psline[linewidth=0.04cm](2.0608783,0.124698274)(2.2846725,0.03800268)
\psline[linewidth=0.04cm](2.0175304,0.01280114)(2.2413247,-0.07389445)
\psline[linewidth=0.04cm](1.966958,-0.11774552)(2.1907523,-0.20444112)
\psline[linewidth=0.04cm](1.9163855,-0.24829218)(2.1401799,-0.33498776)
\psline[linewidth=0.04cm](3.22,0.66)(3.22,-0.54)
\psdots[dotsize=0.24](3.24,0.7)
\pscircle[linewidth=0.03,dimen=outer](3.84,-1.7){0.26}
\psline[linewidth=0.04cm](3.8424585,-1.7610121)(4.4172077,-0.43871742)
\psline[linewidth=0.04cm](3.8424585,-1.7610121)(3.2372172,-0.443407)
\psdots[dotsize=0.24](3.22,-0.5) \psdots[dotsize=0.24](4.4,-0.5)
\psline[linewidth=0.04cm](4.17982,-0.64501226)(4.4036145,-0.7317078)
\psline[linewidth=0.04cm](4.1292477,-0.7755589)(4.353042,-0.8622545)
\psline[linewidth=0.04cm](4.0714507,-0.9247551)(4.2952447,-1.0114506)
\psline[linewidth=0.04cm](4.0208783,-1.0553018)(4.2446723,-1.1419973)
\psline[linewidth=0.04cm](3.9775305,-1.1671989)(4.201325,-1.2538944)
\psline[linewidth=0.04cm](3.926958,-1.2977455)(4.150752,-1.3844411)
\psline[linewidth=0.04cm](3.8763857,-1.4282922)(4.1001797,-1.5149877)
\psline[linewidth=0.04cm](8.282458,-1.781012)(8.857208,-0.4587174)
\psline[linewidth=0.04cm](8.282458,-1.781012)(7.677217,-0.463407)
\psdots[dotsize=0.24](7.66,-0.52) \psdots[dotsize=0.24](8.84,-0.52)
\psline[linewidth=0.08cm](8.12,-1.56)(8.4,-1.86)
\psline[linewidth=0.08cm](8.4,-1.56)(8.12,-1.86)
\pscircle[linewidth=0.03,dimen=outer](8.82,-0.54){0.26}
\psline[linewidth=0.04cm](8.7,0.5)(8.94,0.5)
\psline[linewidth=0.04cm](8.7,0.36)(8.94,0.36)
\psline[linewidth=0.04cm](8.7,0.2)(8.94,0.2)
\psline[linewidth=0.04cm](8.7,0.06)(8.94,0.06)
\psline[linewidth=0.04cm](8.7,-0.06)(8.94,-0.06)
\psline[linewidth=0.04cm](8.7,-0.2)(8.94,-0.2)
\psline[linewidth=0.04cm](8.7,-0.34)(8.94,-0.34)
\psline[linewidth=0.04cm](8.82,0.62)(8.82,-0.58)
\psdots[dotsize=0.24](8.84,0.66)
\psline[linewidth=0.08cm](10.74,-1.58)(11.02,-1.88)
\psline[linewidth=0.08cm](11.02,-1.58)(10.74,-1.88)
\pscircle[linewidth=0.03,dimen=outer](10.88,-1.72){0.26}
\psline[linewidth=0.04cm](10.863444,-1.7736831)(11.88,-0.48)
\psdots[dotsize=0.24,dotangle=-14.541179](11.879736,-0.53304964)
\psline[linewidth=0.04cm](11.470201,-0.77813494)(11.665058,-0.91824275)
\psline[linewidth=0.04cm](11.388471,-0.8918024)(11.583329,-1.0319102)
\psline[linewidth=0.04cm](11.295065,-1.0217079)(11.4899235,-1.1618158)
\psline[linewidth=0.04cm](11.213336,-1.1353754)(11.408195,-1.2754831)
\psline[linewidth=0.04cm](11.143282,-1.2328045)(11.3381405,-1.3729124)
\psline[linewidth=0.04cm](11.061552,-1.3464719)(11.256411,-1.4865798)
\psline[linewidth=0.04cm](10.959823,-1.4401394)(11.154681,-1.5802472)
\psline[linewidth=0.04cm](10.902458,-1.781012)(9.92,-0.56)
\psdots[dotsize=0.24](9.88,-0.54)
\psline[linewidth=0.04cm](11.6502,-0.59813493)(11.845058,-0.73824275)
\psline[linewidth=0.04cm](11.568471,-0.71180236)(11.763329,-0.8519102)
\psline[linewidth=0.04cm](10.899823,-1.5401393)(11.094681,-1.6802472)
\end{pspicture}
}\\\vskip2mm\footnotesize{Figure 1. Trees with the short lived
mutation}
\end{center}
\end{figure}

We note that the total number of vertices is $(n+1){2n-2\choose
n-1}$ for $n\ge1$ and the proportion of new type vertices ${2\over
(n+1){2n-2\choose n-1}}$ approaches $0$ when $n$ gets larger as is
obvious. However, if the horse should resume an active social life
the generating function for the trees becomes $zBC$.

For the second example we look at \textit{toggle trees}. Once the
mutator has a child of the new type, all the later descendants are
also new type. In other words, the first child of the mutator of the
new type plays a role of a `toggle' that divides all children of the
mutator into two groups, those on the left are normal, those on the
right are new type.

\begin{theorem} The number of toggle trees with $n$ edges is ${2n+1\choose n}$. In particular,
the proportion of new type vertices is asymptotically
${1\over2}\sqrt{{\pi\over n}}$.
\end{theorem}

\noindent\textbf{Proof.} The number of toggle trees with a mutator
at the root has the generating function $C^2$, where each of the
first family and the new family contributes a $C$.

If we allow a mutator to be anywhere, applying the uplift principle
gives
\begin{eqnarray*}
T_M={B\over C}\cdot C^2=BC=\sum_{n\ge0}{2n+1\choose
n}z^n=1+3z+10z^2+35z^3+\cdots\quad (A001700).
\end{eqnarray*}

Then the generating function for vertices is
\begin{eqnarray*}
(zT_M)^\prime=\left({B-1\over2}\right)^\prime=B^3=\sum_{n\ge0}(2n+1){2n\choose
n}z^n.
\end{eqnarray*}
To count new type vertices we see that if the mutator is at the
root, we have $C\cdot(zC)^\prime=CB$ possibilities with $C$ for the
pretoggle subtree and $(zC)^\prime=B$ counting the new type
vertices. Multiplying by $L=B/C$ allows the mutator to be anywhere
and our generating function is $V_N={B\over C}\cdot
CB=B^2=1+4z+16z^2+64z^3+\cdots$.

To estimate the proportion of new type vertices, we apply the
Stirling's approximation and get
\begin{eqnarray}\label{e:toggp}
{[z^n]B^2\over [z^n]B^3}={4^n\over (2n+1){2n\choose
n}}\sim{\sqrt{\pi n}{2n\choose n}\over(2n+1){2n\choose
n}}\sim{1\over2}\sqrt{{\pi\over n}}.
\end{eqnarray}
\hfill{\rule{2mm}{2mm}}\bigskip

Figure 2 illustrates the ten toggle trees with $30$ vertices of
which $16$ are of the new type.
\begin{figure}[h]
\begin{center}
\scalebox{0.6} {
\begin{pspicture}(0,-1.46)(14.5,1.46)
\psline[linewidth=0.04cm](0.26,1.18)(0.26,-1.22)
\psdots[dotsize=0.24](0.26,1.2)
\psline[linewidth=0.08cm](0.12,-1.04)(0.4,-1.34)
\psline[linewidth=0.08cm](0.4,-1.04)(0.12,-1.34)
\pscircle[linewidth=0.03,dimen=outer](0.26,-1.2){0.26}
\psdots[dotsize=0.24](0.26,0.0)
\psline[linewidth=0.04cm](1.16,1.02)(1.4,1.02)
\psline[linewidth=0.04cm](1.16,0.88)(1.4,0.88)
\psline[linewidth=0.04cm](1.16,0.72)(1.4,0.72)
\psline[linewidth=0.04cm](1.16,0.58)(1.4,0.58)
\psline[linewidth=0.04cm](1.16,0.46)(1.4,0.46)
\psline[linewidth=0.04cm](1.16,0.32)(1.4,0.32)
\psline[linewidth=0.04cm](1.16,0.18)(1.4,0.18)
\psline[linewidth=0.04cm](1.28,1.18)(1.28,-1.22)
\psdots[dotsize=0.24](1.28,1.2)
\psline[linewidth=0.08cm](1.14,-1.04)(1.42,-1.34)
\psline[linewidth=0.08cm](1.42,-1.04)(1.14,-1.34)
\pscircle[linewidth=0.03,dimen=outer](1.28,-1.2){0.26}
\psdots[dotsize=0.24](1.28,0.0)
\psline[linewidth=0.04cm](1.16,-0.18)(1.4,-0.18)
\psline[linewidth=0.04cm](1.16,-0.32)(1.4,-0.32)
\psline[linewidth=0.04cm](1.16,-0.48)(1.4,-0.48)
\psline[linewidth=0.04cm](1.16,-0.62)(1.4,-0.62)
\psline[linewidth=0.04cm](1.16,-0.74)(1.4,-0.74)
\psline[linewidth=0.04cm](1.16,-0.88)(1.4,-0.88)
\psline[linewidth=0.04cm](1.16,-1.02)(1.4,-1.02)
\psline[linewidth=0.04cm](3.26,1.18)(3.26,-1.22)
\psdots[dotsize=0.24](3.26,1.2)
\psline[linewidth=0.08cm](3.12,-1.04)(3.4,-1.34)
\psline[linewidth=0.08cm](3.4,-1.04)(3.12,-1.34)
\pscircle[linewidth=0.03,dimen=outer](3.26,-0.02){0.26}
\psdots[dotsize=0.24](3.26,0.0)
\psline[linewidth=0.04cm](3.14,1.02)(3.38,1.02)
\psline[linewidth=0.04cm](3.14,0.88)(3.38,0.88)
\psline[linewidth=0.04cm](3.14,0.72)(3.38,0.72)
\psline[linewidth=0.04cm](3.14,0.58)(3.38,0.58)
\psline[linewidth=0.04cm](3.14,0.46)(3.38,0.46)
\psline[linewidth=0.04cm](3.14,0.32)(3.38,0.32)
\psline[linewidth=0.04cm](3.14,0.18)(3.38,0.18)
\psline[linewidth=0.04cm](2.28,1.18)(2.28,-1.22)
\psdots[dotsize=0.24](2.28,1.2)
\psline[linewidth=0.08cm](2.14,-1.04)(2.42,-1.34)
\psline[linewidth=0.08cm](2.42,-1.04)(2.14,-1.34)
\pscircle[linewidth=0.03,dimen=outer](2.28,0.0){0.26}
\psdots[dotsize=0.24](2.28,0.0)
\psline[linewidth=0.04cm](4.24,1.18)(4.24,-1.22)
\psdots[dotsize=0.24](4.24,1.2)
\psline[linewidth=0.08cm](4.1,-1.04)(4.38,-1.34)
\psline[linewidth=0.08cm](4.38,-1.04)(4.1,-1.34)
\pscircle[linewidth=0.03,dimen=outer](4.24,1.2){0.26}
\psdots[dotsize=0.24](4.24,0.0)
\psline[linewidth=0.04cm](5.7224584,-1.2410121)(6.297208,0.08128259)
\psline[linewidth=0.04cm](5.7224584,-1.2410121)(5.117217,0.076593)
\psdots[dotsize=0.24](5.1,0.02) \psdots[dotsize=0.24](6.28,0.02)
\psline[linewidth=0.08cm](5.56,-1.02)(5.84,-1.32)
\psline[linewidth=0.08cm](5.84,-1.02)(5.56,-1.32)
\pscircle[linewidth=0.03,dimen=outer](5.7,-1.2){0.26}
\psline[linewidth=0.08cm](7.56,-1.02)(7.84,-1.32)
\psline[linewidth=0.08cm](7.84,-1.02)(7.56,-1.32)
\pscircle[linewidth=0.03,dimen=outer](7.7,-1.18){0.26}
\psline[linewidth=0.04cm](7.7024584,-1.2410121)(8.277207,0.08128259)
\psline[linewidth=0.04cm](7.7024584,-1.2410121)(7.097217,0.076593)
\psdots[dotsize=0.24](7.08,0.02) \psdots[dotsize=0.24](8.26,0.02)
\psline[linewidth=0.04cm](8.03982,-0.12501223)(8.263615,-0.21170783)
\psline[linewidth=0.04cm](7.989248,-0.25555888)(8.213042,-0.3422545)
\psline[linewidth=0.04cm](7.931451,-0.40475506)(8.155245,-0.49145067)
\psline[linewidth=0.04cm](7.8808784,-0.53530174)(8.104672,-0.6219973)
\psline[linewidth=0.04cm](7.8375306,-0.64719886)(8.061325,-0.73389447)
\psline[linewidth=0.04cm](7.786958,-0.77774554)(8.010753,-0.8644411)
\psline[linewidth=0.04cm](7.7363853,-0.9082922)(7.96018,-0.9949878)
\psline[linewidth=0.08cm](9.56,-1.0)(9.84,-1.3)
\psline[linewidth=0.08cm](9.84,-1.0)(9.56,-1.3)
\pscircle[linewidth=0.03,dimen=outer](9.7,-1.16){0.26}
\psline[linewidth=0.04cm](9.702458,-1.221012)(10.277207,0.10128259)
\psline[linewidth=0.04cm](9.702458,-1.221012)(9.097218,0.096593)
\psdots[dotsize=0.24](9.08,0.04) \psdots[dotsize=0.24](10.26,0.04)
\psline[linewidth=0.04cm](10.03982,-0.10501224)(10.263615,-0.19170783)
\psline[linewidth=0.04cm](9.989247,-0.2355589)(10.213042,-0.32225448)
\psline[linewidth=0.04cm](9.931451,-0.38475507)(10.155245,-0.47145066)
\psline[linewidth=0.04cm](9.880878,-0.5153017)(10.104672,-0.6019973)
\psline[linewidth=0.04cm](9.83753,-0.6271989)(10.061325,-0.7138944)
\psline[linewidth=0.04cm](9.786958,-0.7577455)(10.010753,-0.8444411)
\psline[linewidth=0.04cm](9.736385,-0.8882922)(9.96018,-0.97498775)
\psline[linewidth=0.04cm](9.13694,-0.24731822)(9.355423,-0.14799333)
\psline[linewidth=0.04cm](9.19488,-0.3747663)(9.4133625,-0.2754414)
\psline[linewidth=0.04cm](9.261096,-0.52042127)(9.479579,-0.42109635)
\psline[linewidth=0.04cm](9.319036,-0.6478693)(9.5375185,-0.5485444)
\psline[linewidth=0.04cm](9.368698,-0.75711054)(9.587181,-0.65778565)
\psline[linewidth=0.04cm](9.426638,-0.8845586)(9.645121,-0.7852337)
\psline[linewidth=0.04cm](9.484577,-1.0120066)(9.70306,-0.91268176)
\psline[linewidth=0.04cm](11.702458,-1.221012)(12.277207,0.10128259)
\psline[linewidth=0.04cm](11.702458,-1.221012)(11.097218,0.096593)
\psdots[dotsize=0.24](11.08,0.04) \psdots[dotsize=0.24](12.26,0.04)
\psline[linewidth=0.08cm](11.54,-1.0)(11.82,-1.3)
\psline[linewidth=0.08cm](11.82,-1.0)(11.54,-1.3)
\pscircle[linewidth=0.03,dimen=outer](11.06,0.04){0.26}
\psline[linewidth=0.04cm](13.702458,-1.221012)(14.277207,0.10128259)
\psline[linewidth=0.04cm](13.702458,-1.221012)(13.097218,0.096593)
\psdots[dotsize=0.24](13.08,0.04) \psdots[dotsize=0.24](14.26,0.04)
\psline[linewidth=0.08cm](13.54,-1.0)(13.82,-1.3)
\psline[linewidth=0.08cm](13.82,-1.0)(13.54,-1.3)
\pscircle[linewidth=0.03,dimen=outer](14.24,0.02){0.26}
\end{pspicture}
}\\\vskip2mm\footnotesize{Figure 2. Toggle trees with $16$ new type
vertices.}
\end{center}
\end{figure}
The result \eref{e:toggp} is reasonable since a mutator high up will
usually have few descendants. If, by way of contrast, we specify
that the mutator be at height $1$ then the number of toggle trees is
counted by
\begin{eqnarray*}
zC^4=\sum_{n\ge1}{4\over n+3}{2n+1\choose
n-1}z^n=z+4z^2+14z^3+48z^4+165z^5+\cdots.
\end{eqnarray*}

The number of vertices has the generating function
$(z^2C^4)'=2zC^4B=2\sum_{n\ge1}{2n+2\choose
n-1}z^n=2z+12z^2+56z^3+240z^4+990z^5+\cdots$. The generating
function for vertices of the new type is
\begin{eqnarray*}
zC^3B=\sum_{n\ge1}{2n+1\choose n-1}z^n=z+5z^2+21z^3+84z^4+\cdots.
\end{eqnarray*}
For instance, there are $4\cdot14=56$ vertices of which $21$
vertices are of the new type, see Figure 3.
\begin{figure}[h]
\begin{center}
\scalebox{0.5} {
\begin{pspicture}(0,-1.98)(27.1,1.96)
\psline[linewidth=0.04cm](0.26,1.8)(0.26,-1.8)
\psline[linewidth=0.04cm](8.882459,-1.841012)(9.457208,-0.5187174)
\psline[linewidth=0.04cm](8.882459,-1.841012)(8.277217,-0.523407)
\psdots[dotsize=0.24](0.26,1.82) \psdots[dotsize=0.24](0.26,0.62)
\psline[linewidth=0.08cm](0.12,-1.62)(0.4,-1.92)
\psline[linewidth=0.08cm](0.4,-1.62)(0.12,-1.92)
\psline[linewidth=0.04cm](2.88,-0.6)(2.88,-1.8)
\psdots[dotsize=0.24](2.88,-0.6)
\psline[linewidth=0.08cm](2.72,-1.62)(3.0,-1.92)
\psline[linewidth=0.08cm](3.0,-1.62)(2.72,-1.92)
\psdots[dotsize=0.24](8.26,-0.58) \psdots[dotsize=0.24](9.44,-0.58)
\psline[linewidth=0.08cm](8.72,-1.62)(9.0,-1.92)
\psline[linewidth=0.08cm](9.0,-1.62)(8.72,-1.92)
\psline[linewidth=0.04cm](21.26,-1.7610121)(21.26,-0.56)
\psdots[dotsize=0.24](21.26,-0.56)
\pscircle[linewidth=0.03,dimen=outer](0.26,-0.58){0.26}
\pscircle[linewidth=0.03,dimen=outer](2.86,-0.6){0.26}
\psdots[dotsize=0.24](0.26,-0.58)
\pscircle[linewidth=0.03,dimen=outer](8.24,-0.58){0.26}
\psline[linewidth=0.04cm](8.24,0.58)(8.24,-0.62)
\psdots[dotsize=0.24](8.26,0.62)
\psline[linewidth=0.04cm](2.9024584,-0.66101205)(3.477208,0.6612826)
\psline[linewidth=0.04cm](2.9024584,-0.66101205)(2.2972171,0.656593)
\psdots[dotsize=0.24](2.28,0.6) \psdots[dotsize=0.24](3.46,0.6)
\psline[linewidth=0.04cm](14.902458,-1.821012)(15.477208,-0.4987174)
\psline[linewidth=0.04cm](14.902458,-1.821012)(14.297217,-0.503407)
\psdots[dotsize=0.24](14.28,-0.56)
\psdots[dotsize=0.24](15.46,-0.56)
\psline[linewidth=0.08cm](14.74,-1.6)(15.02,-1.9)
\psline[linewidth=0.08cm](15.02,-1.6)(14.74,-1.9)
\pscircle[linewidth=0.03,dimen=outer](14.26,-0.54){0.26}
\psline[linewidth=0.04cm](15.44,0.58)(15.44,-0.62)
\psdots[dotsize=0.24](15.46,0.62)
\psline[linewidth=0.08cm](21.12,-1.62)(21.4,-1.92)
\psline[linewidth=0.08cm](21.4,-1.62)(21.12,-1.92)
\pscircle[linewidth=0.03,dimen=outer](20.44,-0.58){0.26}
\psline[linewidth=0.04cm](21.243444,-1.813683)(22.04,-0.58)
\psdots[dotsize=0.24,dotangle=-14.541179](22.059736,-0.57304966)
\psline[linewidth=0.04cm](21.28246,-1.821012)(20.48,-0.64)
\psdots[dotsize=0.24](20.46,-0.58)
\psline[linewidth=0.04cm](1.26,1.8)(1.26,-1.8)
\psdots[dotsize=0.24](1.26,1.82) \psdots[dotsize=0.24](1.26,0.62)
\psline[linewidth=0.08cm](1.12,-1.62)(1.4,-1.92)
\psline[linewidth=0.08cm](1.4,-1.62)(1.12,-1.92)
\pscircle[linewidth=0.03,dimen=outer](1.26,-0.58){0.26}
\psdots[dotsize=0.24](1.26,-0.58)
\psline[linewidth=0.04cm](1.14,1.66)(1.38,1.66)
\psline[linewidth=0.04cm](1.14,1.52)(1.38,1.52)
\psline[linewidth=0.04cm](1.14,1.36)(1.38,1.36)
\psline[linewidth=0.04cm](1.14,1.22)(1.38,1.22)
\psline[linewidth=0.04cm](1.14,1.1)(1.38,1.1)
\psline[linewidth=0.04cm](1.14,0.96)(1.38,0.96)
\psline[linewidth=0.04cm](1.14,0.82)(1.38,0.82)
\psline[linewidth=0.04cm](1.14,0.46)(1.38,0.46)
\psline[linewidth=0.04cm](1.14,0.32)(1.38,0.32)
\psline[linewidth=0.04cm](1.14,0.16)(1.38,0.16)
\psline[linewidth=0.04cm](1.14,0.02)(1.38,0.02)
\psline[linewidth=0.04cm](1.14,-0.1)(1.38,-0.1)
\psline[linewidth=0.04cm](1.14,-0.24)(1.38,-0.24)
\psline[linewidth=0.04cm](1.14,-0.38)(1.38,-0.38)
\psline[linewidth=0.04cm](4.88,-0.62)(4.88,-1.82)
\psdots[dotsize=0.24](4.88,-0.62)
\psline[linewidth=0.08cm](4.72,-1.64)(5.0,-1.94)
\psline[linewidth=0.08cm](5.0,-1.64)(4.72,-1.94)
\pscircle[linewidth=0.03,dimen=outer](4.86,-0.62){0.26}
\psline[linewidth=0.04cm](4.902458,-0.68101203)(5.4772077,0.6412826)
\psline[linewidth=0.04cm](4.902458,-0.68101203)(4.2972174,0.636593)
\psdots[dotsize=0.24](4.28,0.58) \psdots[dotsize=0.24](5.46,0.58)
\psline[linewidth=0.04cm](5.23982,0.43498775)(5.4636145,0.34829217)
\psline[linewidth=0.04cm](5.1892476,0.3044411)(5.413042,0.21774551)
\psline[linewidth=0.04cm](5.1314507,0.15524493)(5.355245,0.068549335)
\psline[linewidth=0.04cm](5.0808783,0.024698274)(5.3046727,-0.061997317)
\psline[linewidth=0.04cm](5.0375304,-0.08719886)(5.261325,-0.17389445)
\psline[linewidth=0.04cm](4.986958,-0.21774551)(5.2107525,-0.3044411)
\psline[linewidth=0.04cm](4.9363856,-0.34829217)(5.1601796,-0.43498775)
\psline[linewidth=0.04cm](6.88,-0.6)(6.88,-1.8)
\psdots[dotsize=0.24](6.88,-0.6)
\psline[linewidth=0.08cm](6.72,-1.62)(7.0,-1.92)
\psline[linewidth=0.08cm](7.0,-1.62)(6.72,-1.92)
\pscircle[linewidth=0.03,dimen=outer](6.86,-0.6){0.26}
\psline[linewidth=0.04cm](6.902458,-0.66101205)(7.4772077,0.6612826)
\psline[linewidth=0.04cm](6.902458,-0.66101205)(6.2972174,0.656593)
\psdots[dotsize=0.24](6.28,0.6) \psdots[dotsize=0.24](7.46,0.6)
\psline[linewidth=0.04cm](7.23982,0.45498776)(7.4636145,0.36829218)
\psline[linewidth=0.04cm](7.1892476,0.3244411)(7.413042,0.23774552)
\psline[linewidth=0.04cm](7.1314507,0.17524493)(7.355245,0.08854934)
\psline[linewidth=0.04cm](7.0808783,0.044698272)(7.3046727,-0.041997317)
\psline[linewidth=0.04cm](7.0375304,-0.06719886)(7.261325,-0.15389445)
\psline[linewidth=0.04cm](6.986958,-0.19774552)(7.2107525,-0.2844411)
\psline[linewidth=0.04cm](6.9363856,-0.32829216)(7.1601796,-0.41498777)
\psline[linewidth=0.04cm](6.3369403,0.31268176)(6.555423,0.41200668)
\psline[linewidth=0.04cm](6.39488,0.1852337)(6.6133623,0.2845586)
\psline[linewidth=0.04cm](6.4610963,0.03957876)(6.679579,0.13890365)
\psline[linewidth=0.04cm](6.519036,-0.08786932)(6.7375183,0.01145558)
\psline[linewidth=0.04cm](6.5686984,-0.19711052)(6.787181,-0.09778563)
\psline[linewidth=0.04cm](6.626638,-0.3245586)(6.8451204,-0.2252337)
\psline[linewidth=0.04cm](6.6845775,-0.45200667)(6.90306,-0.3526818)
\psline[linewidth=0.04cm](10.882459,-1.821012)(11.457208,-0.4987174)
\psline[linewidth=0.04cm](10.882459,-1.821012)(10.277217,-0.503407)
\psdots[dotsize=0.24](10.26,-0.56)
\psdots[dotsize=0.24](11.44,-0.56)
\psline[linewidth=0.08cm](10.72,-1.6)(11.0,-1.9)
\psline[linewidth=0.08cm](11.0,-1.6)(10.72,-1.9)
\pscircle[linewidth=0.03,dimen=outer](10.24,-0.56){0.26}
\psline[linewidth=0.04cm](10.12,0.48)(10.36,0.48)
\psline[linewidth=0.04cm](10.12,0.34)(10.36,0.34)
\psline[linewidth=0.04cm](10.12,0.18)(10.36,0.18)
\psline[linewidth=0.04cm](10.12,0.04)(10.36,0.04)
\psline[linewidth=0.04cm](10.12,-0.08)(10.36,-0.08)
\psline[linewidth=0.04cm](10.12,-0.22)(10.36,-0.22)
\psline[linewidth=0.04cm](10.12,-0.36)(10.36,-0.36)
\psline[linewidth=0.04cm](10.24,0.6)(10.24,-0.6)
\psdots[dotsize=0.24](10.26,0.64)
\psline[linewidth=0.04cm](12.882459,-1.821012)(13.457208,-0.4987174)
\psline[linewidth=0.04cm](12.882459,-1.821012)(12.277217,-0.503407)
\psdots[dotsize=0.24](12.26,-0.56)
\psdots[dotsize=0.24](13.44,-0.56)
\psline[linewidth=0.08cm](12.72,-1.6)(13.0,-1.9)
\psline[linewidth=0.08cm](13.0,-1.6)(12.72,-1.9)
\pscircle[linewidth=0.03,dimen=outer](13.42,-0.56){0.26}
\psline[linewidth=0.04cm](12.24,0.6)(12.24,-0.6)
\psdots[dotsize=0.24](12.26,0.64)
\psline[linewidth=0.04cm](16.86246,-1.821012)(17.437208,-0.4987174)
\psline[linewidth=0.04cm](16.86246,-1.821012)(16.257217,-0.503407)
\psdots[dotsize=0.24](16.24,-0.56)
\psdots[dotsize=0.24](17.42,-0.56)
\psline[linewidth=0.08cm](16.7,-1.6)(16.98,-1.9)
\psline[linewidth=0.08cm](16.98,-1.6)(16.7,-1.9)
\pscircle[linewidth=0.03,dimen=outer](17.4,-0.58){0.26}
\psline[linewidth=0.04cm](17.4,0.58)(17.4,-0.62)
\psdots[dotsize=0.24](17.42,0.62)
\psline[linewidth=0.04cm](18.882458,-1.821012)(19.457209,-0.4987174)
\psline[linewidth=0.04cm](18.882458,-1.821012)(18.277218,-0.503407)
\psdots[dotsize=0.24](18.26,-0.56)
\psdots[dotsize=0.24](19.44,-0.56)
\psline[linewidth=0.08cm](18.72,-1.6)(19.0,-1.9)
\psline[linewidth=0.08cm](19.0,-1.6)(18.72,-1.9)
\pscircle[linewidth=0.03,dimen=outer](19.42,-0.58){0.26}
\psline[linewidth=0.04cm](19.3,0.46)(19.54,0.46)
\psline[linewidth=0.04cm](19.3,0.32)(19.54,0.32)
\psline[linewidth=0.04cm](19.3,0.16)(19.54,0.16)
\psline[linewidth=0.04cm](19.3,0.02)(19.54,0.02)
\psline[linewidth=0.04cm](19.3,-0.1)(19.54,-0.1)
\psline[linewidth=0.04cm](19.3,-0.24)(19.54,-0.24)
\psline[linewidth=0.04cm](19.3,-0.38)(19.54,-0.38)
\psline[linewidth=0.04cm](19.42,0.58)(19.42,-0.62)
\psdots[dotsize=0.24](19.44,0.62)
\psline[linewidth=0.04cm](23.68,-1.7410121)(23.68,-0.54)
\psdots[dotsize=0.24](23.68,-0.54)
\psline[linewidth=0.08cm](23.54,-1.6)(23.82,-1.9)
\psline[linewidth=0.08cm](23.82,-1.6)(23.54,-1.9)
\pscircle[linewidth=0.03,dimen=outer](23.66,-0.52){0.26}
\psline[linewidth=0.04cm](23.663445,-1.793683)(24.46,-0.56)
\psdots[dotsize=0.24,dotangle=-14.541179](24.479736,-0.5530496)
\psline[linewidth=0.04cm](23.70246,-1.801012)(22.9,-0.62)
\psdots[dotsize=0.24](22.88,-0.56)
\psline[linewidth=0.04cm](26.06,-1.7410121)(26.06,-0.54)
\psdots[dotsize=0.24](26.06,-0.54)
\psline[linewidth=0.08cm](25.92,-1.6)(26.2,-1.9)
\psline[linewidth=0.08cm](26.2,-1.6)(25.92,-1.9)
\pscircle[linewidth=0.03,dimen=outer](26.84,-0.54){0.26}
\psline[linewidth=0.04cm](26.043444,-1.793683)(26.84,-0.56)
\psdots[dotsize=0.24,dotangle=-14.541179](26.859735,-0.5530496)
\psline[linewidth=0.04cm](26.082458,-1.801012)(25.28,-0.62)
\psdots[dotsize=0.24](25.26,-0.56)
\end{pspicture}
}\\\vskip2mm\footnotesize{Figure 3. Toggle trees with a mutator at
height $1$.}
\end{center}
\end{figure}

\noindent We note that the proportion of new type vertices in this
height $1$ case is
\begin{eqnarray*}
{{2n+1\choose n-1}\over {4n+4\over n+3}{2n+1\choose n-1}}={n+3\over
4n+4}\rightarrow{1\over4}.
\end{eqnarray*}

In a toggle tree, if a child of the new type of the mutator appears,
all later offspring of the mutator and their descendants are of the
new type. Suppose instead that every child of the mutator and
recursively every child of a new type vertex has a $50\%$ chance of
being new type. We call such trees \textit{embedded new type (ENT)
trees} since the result is a tree with a subtree of new type
vertices. This concept coincides with an autosomal dominant mutation
when we assume that we only keep track of genetic history of a
single family, and every member met a spouse not having a mutant
gene so that appearance of a mutation only depends on a member of
the family. Figure 4 illustrates the $12$ possible trees with $2$
edges.

\begin{center}
\scalebox{0.6} {
\begin{pspicture}(0,-1.46)(17.64,1.46)
\psline[linewidth=0.04cm](0.26,1.18)(0.26,-1.22)
\psdots[dotsize=0.24](0.26,1.2) \psdots[dotsize=0.24](0.26,0.0)
\psline[linewidth=0.08cm](0.12,-1.04)(0.4,-1.34)
\psline[linewidth=0.08cm](0.4,-1.04)(0.12,-1.34)
\psline[linewidth=0.08cm](6.72,-1.04)(7.0,-1.34)
\psline[linewidth=0.08cm](7.0,-1.04)(6.72,-1.34)
\pscircle[linewidth=0.03,dimen=outer](0.26,-1.2){0.26}
\pscircle[linewidth=0.03,dimen=outer](6.86,-1.18){0.26}
\psline[linewidth=0.04cm](6.902458,-1.2410121)(7.4772077,0.08128259)
\psline[linewidth=0.04cm](6.902458,-1.2410121)(6.2972174,0.076593)
\psdots[dotsize=0.24](6.28,0.02) \psdots[dotsize=0.24](7.46,0.02)
\psline[linewidth=0.04cm](14.882459,-1.2410121)(15.457208,0.08128259)
\psline[linewidth=0.04cm](14.882459,-1.2410121)(14.277217,0.076593)
\psdots[dotsize=0.24](14.26,0.02) \psdots[dotsize=0.24](15.44,0.02)
\psline[linewidth=0.08cm](14.72,-1.02)(15.0,-1.32)
\psline[linewidth=0.08cm](15.0,-1.02)(14.72,-1.32)
\pscircle[linewidth=0.03,dimen=outer](14.24,0.04){0.26}
\psline[linewidth=0.04cm](1.26,1.18)(1.26,-1.22)
\psdots[dotsize=0.24](1.26,1.2) \psdots[dotsize=0.24](1.26,0.0)
\psline[linewidth=0.08cm](1.12,-1.04)(1.4,-1.34)
\psline[linewidth=0.08cm](1.4,-1.04)(1.12,-1.34)
\pscircle[linewidth=0.03,dimen=outer](1.26,-1.2){0.26}
\psline[linewidth=0.04cm](1.14,-0.16)(1.38,-0.16)
\psline[linewidth=0.04cm](1.14,-0.3)(1.38,-0.3)
\psline[linewidth=0.04cm](1.14,-0.46)(1.38,-0.46)
\psline[linewidth=0.04cm](1.14,-0.6)(1.38,-0.6)
\psline[linewidth=0.04cm](1.14,-0.72)(1.38,-0.72)
\psline[linewidth=0.04cm](1.14,-0.86)(1.38,-0.86)
\psline[linewidth=0.04cm](1.14,-1.0)(1.38,-1.0)
\psline[linewidth=0.08cm](10.72,-1.04)(11.0,-1.34)
\psline[linewidth=0.08cm](11.0,-1.04)(10.72,-1.34)
\pscircle[linewidth=0.03,dimen=outer](10.86,-1.18){0.26}
\psline[linewidth=0.04cm](10.902458,-1.2410121)(11.477208,0.08128259)
\psline[linewidth=0.04cm](10.902458,-1.2410121)(10.297217,0.076593)
\psdots[dotsize=0.24](10.28,0.02) \psdots[dotsize=0.24](11.46,0.02)
\psline[linewidth=0.04cm](11.2398205,-0.12501223)(11.463614,-0.21170783)
\psline[linewidth=0.04cm](11.189248,-0.25555888)(11.413042,-0.3422545)
\psline[linewidth=0.04cm](11.131451,-0.40475506)(11.355245,-0.49145067)
\psline[linewidth=0.04cm](11.080878,-0.53530174)(11.304672,-0.6219973)
\psline[linewidth=0.04cm](11.037531,-0.64719886)(11.261325,-0.73389447)
\psline[linewidth=0.04cm](10.986958,-0.77774554)(11.2107525,-0.8644411)
\psline[linewidth=0.04cm](10.936385,-0.9082922)(11.16018,-0.9949878)
\psline[linewidth=0.08cm](12.72,-1.02)(13.0,-1.32)
\psline[linewidth=0.08cm](13.0,-1.02)(12.72,-1.32)
\pscircle[linewidth=0.03,dimen=outer](12.86,-1.16){0.26}
\psline[linewidth=0.04cm](12.902458,-1.221012)(13.477208,0.10128259)
\psline[linewidth=0.04cm](12.902458,-1.221012)(12.297217,0.096593)
\psdots[dotsize=0.24](12.28,0.04) \psdots[dotsize=0.24](13.46,0.04)
\psline[linewidth=0.04cm](13.2398205,-0.10501224)(13.463614,-0.19170783)
\psline[linewidth=0.04cm](13.189248,-0.2355589)(13.413042,-0.32225448)
\psline[linewidth=0.04cm](13.131451,-0.38475507)(13.355245,-0.47145066)
\psline[linewidth=0.04cm](13.080878,-0.5153017)(13.304672,-0.6019973)
\psline[linewidth=0.04cm](13.037531,-0.6271989)(13.261325,-0.7138944)
\psline[linewidth=0.04cm](12.986958,-0.7577455)(13.2107525,-0.8444411)
\psline[linewidth=0.04cm](12.936385,-0.8882922)(13.16018,-0.97498775)
\psline[linewidth=0.04cm](12.33694,-0.24731822)(12.555423,-0.14799333)
\psline[linewidth=0.04cm](12.394879,-0.3747663)(12.613362,-0.2754414)
\psline[linewidth=0.04cm](12.461097,-0.52042127)(12.679579,-0.42109635)
\psline[linewidth=0.04cm](12.519036,-0.6478693)(12.737518,-0.5485444)
\psline[linewidth=0.04cm](12.568698,-0.75711054)(12.787181,-0.65778565)
\psline[linewidth=0.04cm](12.626637,-0.8845586)(12.84512,-0.7852337)
\psline[linewidth=0.04cm](12.684577,-1.0120066)(12.90306,-0.91268176)
\psline[linewidth=0.04cm](16.842459,-1.2410121)(17.417208,0.08128259)
\psline[linewidth=0.04cm](16.842459,-1.2410121)(16.237217,0.076593)
\psdots[dotsize=0.24](16.22,0.02) \psdots[dotsize=0.24](17.4,0.02)
\psline[linewidth=0.08cm](16.68,-1.02)(16.96,-1.32)
\psline[linewidth=0.08cm](16.96,-1.02)(16.68,-1.32)
\pscircle[linewidth=0.03,dimen=outer](17.38,0.0){0.26}
\psline[linewidth=0.04cm](2.28,1.18)(2.28,-1.22)
\psdots[dotsize=0.24](2.28,1.2) \psdots[dotsize=0.24](2.28,0.0)
\psline[linewidth=0.08cm](2.14,-1.04)(2.42,-1.34)
\psline[linewidth=0.08cm](2.42,-1.04)(2.14,-1.34)
\pscircle[linewidth=0.03,dimen=outer](2.28,-1.2){0.26}
\psline[linewidth=0.04cm](2.16,1.04)(2.4,1.04)
\psline[linewidth=0.04cm](2.16,0.9)(2.4,0.9)
\psline[linewidth=0.04cm](2.16,0.74)(2.4,0.74)
\psline[linewidth=0.04cm](2.16,0.6)(2.4,0.6)
\psline[linewidth=0.04cm](2.16,0.48)(2.4,0.48)
\psline[linewidth=0.04cm](2.16,0.34)(2.4,0.34)
\psline[linewidth=0.04cm](2.16,0.2)(2.4,0.2)
\psline[linewidth=0.04cm](2.16,-0.16)(2.4,-0.16)
\psline[linewidth=0.04cm](2.16,-0.3)(2.4,-0.3)
\psline[linewidth=0.04cm](2.16,-0.46)(2.4,-0.46)
\psline[linewidth=0.04cm](2.16,-0.6)(2.4,-0.6)
\psline[linewidth=0.04cm](2.16,-0.72)(2.4,-0.72)
\psline[linewidth=0.04cm](2.16,-0.86)(2.4,-0.86)
\psline[linewidth=0.04cm](2.16,-1.0)(2.4,-1.0)
\psline[linewidth=0.04cm](3.28,1.18)(3.28,-1.22)
\psdots[dotsize=0.24](3.28,1.2) \psdots[dotsize=0.24](3.28,0.0)
\psline[linewidth=0.08cm](3.14,-1.04)(3.42,-1.34)
\psline[linewidth=0.08cm](3.42,-1.04)(3.14,-1.34)
\pscircle[linewidth=0.03,dimen=outer](3.28,0.0){0.26}
\psline[linewidth=0.04cm](4.28,1.18)(4.28,-1.22)
\psdots[dotsize=0.24](4.28,1.2) \psdots[dotsize=0.24](4.28,0.0)
\psline[linewidth=0.08cm](4.14,-1.04)(4.42,-1.34)
\psline[linewidth=0.08cm](4.42,-1.04)(4.14,-1.34)
\pscircle[linewidth=0.03,dimen=outer](4.28,0.0){0.26}
\psline[linewidth=0.04cm](4.16,1.04)(4.4,1.04)
\psline[linewidth=0.04cm](4.16,0.9)(4.4,0.9)
\psline[linewidth=0.04cm](4.16,0.74)(4.4,0.74)
\psline[linewidth=0.04cm](4.16,0.6)(4.4,0.6)
\psline[linewidth=0.04cm](4.16,0.48)(4.4,0.48)
\psline[linewidth=0.04cm](4.16,0.34)(4.4,0.34)
\psline[linewidth=0.04cm](4.16,0.2)(4.4,0.2)
\psline[linewidth=0.04cm](5.3,1.18)(5.3,-1.22)
\psdots[dotsize=0.24](5.3,1.2) \psdots[dotsize=0.24](5.3,0.0)
\psline[linewidth=0.08cm](5.16,-1.04)(5.44,-1.34)
\psline[linewidth=0.08cm](5.44,-1.04)(5.16,-1.34)
\pscircle[linewidth=0.03,dimen=outer](5.28,1.2){0.26}
\psline[linewidth=0.04cm](12.26,-0.14)(12.48,-0.02)
\psline[linewidth=0.08cm](8.72,-1.04)(9.0,-1.34)
\psline[linewidth=0.08cm](9.0,-1.04)(8.72,-1.34)
\pscircle[linewidth=0.03,dimen=outer](8.86,-1.18){0.26}
\psline[linewidth=0.04cm](8.902458,-1.2410121)(9.477208,0.08128259)
\psline[linewidth=0.04cm](8.902458,-1.2410121)(8.297217,0.076593)
\psdots[dotsize=0.24](8.28,0.02) \psdots[dotsize=0.24](9.46,0.02)
\psline[linewidth=0.04cm](8.33694,-0.26731822)(8.555423,-0.16799332)
\psline[linewidth=0.04cm](8.394879,-0.3947663)(8.613362,-0.2954414)
\psline[linewidth=0.04cm](8.461097,-0.54042125)(8.679579,-0.44109634)
\psline[linewidth=0.04cm](8.519036,-0.6678693)(8.737518,-0.56854445)
\psline[linewidth=0.04cm](8.568698,-0.7771105)(8.787181,-0.67778563)
\psline[linewidth=0.04cm](8.626637,-0.9045586)(8.84512,-0.8052337)
\psline[linewidth=0.04cm](8.684577,-1.0320066)(8.90306,-0.9326818)
\psline[linewidth=0.04cm](8.26,-0.16)(8.48,-0.04)
\end{pspicture}
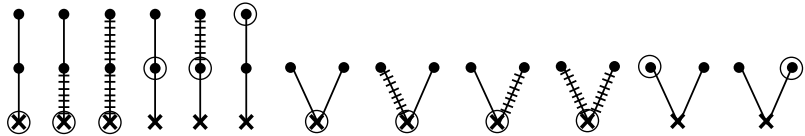
}\\\vskip2mm\footnotesize{Figure 4. The embedded new type trees with
$2$ edges.}
\end{center}

\begin{theorem}
The number of embedded new type trees with $n$ edges is
$\sum_{k=0}^n{1\over k+1}{2n\choose n-k}{2k\choose k}$. In
particular, the proportion of new type vertices is asymptotically
${3\over5}$.
\end{theorem}

\noindent\textbf{Proof.} Let $T_0$ be the generating function for
ENT trees having a mutator at the root. If a mutator has $k$
children, there are $2^k$ possible distributions of the mutation
over the children where each normal child and child of the new type
are the roots of subtrees described by $C$ and $T_0$, respectively.
It then follows that the generating function for ENT trees having a
mutator at the root of updegree $k$ is $z^k(C+T_0)^k$. So $T_0$
satisfies
\begin{eqnarray*}
T_0={1\over1-z(C+T_0)}
\end{eqnarray*}
and solving the functional equation gives
\begin{eqnarray*}
T_0={1-\sqrt{5-4C}\over2zC}=1+2z+7z^2+29z^3+131z^4+\cdots\quad\hbox{(A007852)}.
\end{eqnarray*}
A simple computation shows $T_0=C\cdot(C\circ zC^2)$. By the uplift
principle,
\begin{eqnarray*}
T_M={B\over C}\cdot T_0=B\cdot(C\circ
zC^2)=1+3z+12z^2+52z^3+236z^4+\cdots.
\end{eqnarray*}
This is (A007856) in the OEIS \cite{BSloa} and is also known
\cite{BKlaz} to count the number of subtrees in ordered trees with
$n$ edges. It can be shown that $[z^n]T_M=\sum_{k=0}^nC_k{2n\choose
n-k}$. The singular expansion of $C$ at the dominant singularity
$z={4\over25}$ of $C\circ zC^2$ gives $[z^n]T_M\sim
{5\sqrt{15}\over9}{1\over\sqrt{\pi n^3}}\left({25\over4}\right)^n$,
which implies the asymptotic number of vertices of ENT trees with
$n$ edges is ${5\sqrt{15}\over9}{1\over\sqrt{\pi
n}}\left({25\over4}\right)^n$.

In order to count the vertices of the new type, first consider ENT
trees with a mutator at the root. Let $\tilde{V}_N$ be the
generating function for such ENT trees where we have marked one of
the new type vertices. Suppose the root degree is $k$ of which $j$
are new type. The marked vertex, if not the root itself, is in one
of the $j$ new type subtrees. The other $j-1$ subtrees are
themselves ENT trees with the mutator at the root. Thus
\begin{eqnarray*}
\tilde{V}_N&=&T_0+z\tilde{V}_N+2z^2\tilde{V}_N(C+T_0)+3z^3\tilde{V}_N(C+T_0)^2+\cdots\\
&=&T_0+{z\tilde{V}_N\over(1-z(C+T_0))^2}
\end{eqnarray*}
with $T_0$ counting the case where the marked vertex is the root.
Solving the functional equation yields
\begin{eqnarray*}
V_N&=&L\cdot\tilde{V}_N={B\over C}\cdot
T_0\cdot\left(1-{z\over(1-z(C+T_0))^2}\right)^{-1}=\sqrt{{2-5z+2\sqrt{1-4z}\over(4-25z)(1-4z)}}\\
&=&1+4z+20z^2+106z^3+580z^4+\cdots.
\end{eqnarray*}
It follows from $V_N\sim
{\sqrt{15}\over3}\left(1-{25\over4}z\right)^{-1/2}$ that
$[z^n]V_N\sim{\sqrt{15}\over3}{1\over\sqrt{\pi
n}}\left({25\over4}\right)^n$. Thus the proportion of vertices of
the new type is asymptotically equal to
\begin{eqnarray*}
{{\sqrt{15}\over3}{1\over\sqrt{\pi
n}}\left({25\over4}\right)^n\over{5\sqrt{15}\over9}{1\over\sqrt{\pi
n}}\left({25\over4}\right)^n}={3\over5}.
\end{eqnarray*}
\hfill{\rule{2mm}{2mm}}\bigskip

\noindent In fact, we see for example that
\begin{eqnarray*}
{[z^{50}]V_N\over[z^{50}](zT_M)^\prime}={642784246122173091957609761927581466320\over1086960365349
865718238126127455484769220}\fallingdotseq 0.59135942.
\end{eqnarray*}

One instance where ENT trees occur is classical. These were family
names as in England. The mutator, often the root, passes his name to
his male children (these would be new type) whereas the female
children would not carry on the family name, see \cite{BWat}.

In a \textit{right branch new type (RBNT) tree}, the rightmost
branch from the mutator is nontrivial and all the vertices of the
new type constitute this branch including the mutator.

\begin{center}
\scalebox{0.5} {
\begin{pspicture}(0,-3.54)(6.2115827,3.52)
\psline[linewidth=0.04cm](4.2115827,1.32)(4.2115827,0.12)
\psdots[dotsize=0.24](4.2115827,1.34)
\psline[linewidth=0.04cm](4.091583,1.16)(4.3315825,1.16)
\psline[linewidth=0.04cm](4.091583,1.02)(4.3315825,1.02)
\psline[linewidth=0.04cm](4.091583,0.86)(4.3315825,0.86)
\psline[linewidth=0.04cm](4.091583,0.72)(4.3315825,0.72)
\psline[linewidth=0.04cm](4.091583,0.6)(4.3315825,0.6)
\psline[linewidth=0.04cm](4.091583,0.46)(4.3315825,0.46)
\psline[linewidth=0.04cm](4.091583,0.32)(4.3315825,0.32)
\psline[linewidth=0.04cm](1.1515826,-1.02)(1.1515826,-2.22)
\psdots[dotsize=0.24](1.1515826,-1.0)
\psline[linewidth=0.04cm](2.3515825,-3.401012)(2.3515825,-2.22)
\psline[linewidth=0.04cm](2.3715825,-3.36)(1.1687998,-2.16)
\psdots[dotsize=0.24](1.1515826,-2.22)
\psdots[dotsize=0.24](2.3715825,-2.24)
\psline[linewidth=0.08cm](2.2115827,-3.2)(2.4915826,-3.5)
\psline[linewidth=0.08cm](2.4915826,-3.2)(2.2115827,-3.5)
\psline[linewidth=0.04cm](2.394041,-3.341012)(3.5915825,-2.2)
\psdots[dotsize=0.24](3.5915825,-2.22)
\psline[linewidth=0.04cm](3.614041,-2.241012)(4.1887903,-0.9187174)
\psline[linewidth=0.04cm](3.614041,-2.241012)(3.0087998,-0.923407)
\psdots[dotsize=0.24](2.9915826,-0.98)
\psdots[dotsize=0.24](4.1715827,-0.98)
\psline[linewidth=0.04cm](1.7715826,1.34)(1.7715826,0.14)
\psdots[dotsize=0.24](1.7715826,1.36)
\psline[linewidth=0.04cm](2.9715827,-1.041012)(2.9715827,0.14)
\psline[linewidth=0.04cm](2.9915826,-1.0)(1.7887998,0.2)
\psdots[dotsize=0.24](1.7715826,0.14)
\psdots[dotsize=0.24](2.9915826,0.12)
\psline[linewidth=0.04cm](3.014041,-0.98101205)(4.2115827,0.16)
\psdots[dotsize=0.24](4.2115827,0.14)
\psline[linewidth=0.04cm](3.6948164,-0.14001839)(3.8668287,-0.30738547)
\psline[linewidth=0.04cm](3.5971854,-0.24035895)(3.769198,-0.40772605)
\psline[linewidth=0.04cm](3.4856074,-0.35503387)(3.6576197,-0.522401)
\psline[linewidth=0.04cm](3.3879766,-0.45537442)(3.559989,-0.6227415)
\psline[linewidth=0.04cm](3.3042932,-0.5413806)(3.4763055,-0.7087477)
\psline[linewidth=0.04cm](3.2066622,-0.6417212)(3.3786745,-0.80908823)
\psline[linewidth=0.04cm](3.1090314,-0.74206173)(3.2810438,-0.90942883)
\psline[linewidth=0.04cm](3.8153996,-0.060900368)(3.987412,-0.22826746)
\psline[linewidth=0.04cm](3.9115489,0.03888215)(4.0835614,-0.12848493)
\psline[linewidth=0.04cm](4.214041,1.298988)(4.78879,2.6212826)
\psline[linewidth=0.04cm](4.214041,1.298988)(3.6087997,2.616593)
\psdots[dotsize=0.24](3.5915825,2.56)
\psdots[dotsize=0.24](4.7715826,2.56)
\psline[linewidth=0.04cm](4.5514026,2.4149878)(4.775197,2.3282921)
\psline[linewidth=0.04cm](4.50083,2.284441)(4.7246246,2.1977456)
\psline[linewidth=0.04cm](4.443033,2.1352448)(4.6668277,2.0485494)
\psline[linewidth=0.04cm](4.392461,2.0046983)(4.6162553,1.9180027)
\psline[linewidth=0.04cm](4.349113,1.8928012)(4.5729074,1.8061055)
\psline[linewidth=0.04cm](4.2985406,1.7622545)(4.522335,1.6755589)
\psline[linewidth=0.04cm](4.247968,1.6317078)(4.4717627,1.5450122)
\psline[linewidth=0.04cm](3.6485229,2.2726817)(3.8670053,2.3720067)
\psline[linewidth=0.04cm](3.7064624,2.1452336)(3.9249449,2.2445586)
\psline[linewidth=0.04cm](3.7726789,1.9995787)(3.9911613,2.0989037)
\psline[linewidth=0.04cm](3.8306184,1.8721306)(4.049101,1.9714556)
\psline[linewidth=0.04cm](3.880281,1.7628895)(4.0987635,1.8622143)
\psline[linewidth=0.04cm](3.9382205,1.6354414)(4.156703,1.7347662)
\psline[linewidth=0.04cm](3.99616,1.5079933)(4.2146425,1.6073183)
\psline[linewidth=0.04cm](4.214041,0.15898795)(5.3515825,1.3)
\psdots[dotsize=0.24](5.3515825,1.32)
\psline[linewidth=0.04cm](4.8948164,0.9999816)(5.0668287,0.83261454)
\psline[linewidth=0.04cm](4.7971854,0.89964104)(4.9691978,0.73227394)
\psline[linewidth=0.04cm](4.6856074,0.7849661)(4.85762,0.61759907)
\psline[linewidth=0.04cm](4.5879765,0.68462557)(4.759989,0.51725847)
\psline[linewidth=0.04cm](4.504293,0.5986194)(4.6763053,0.4312523)
\psline[linewidth=0.04cm](4.4066625,0.49827883)(4.578675,0.33091176)
\psline[linewidth=0.04cm](4.3090315,0.39793828)(4.481044,0.2305712)
\psline[linewidth=0.04cm](5.0153995,1.0790997)(5.187412,0.91173255)
\psline[linewidth=0.04cm](5.111549,1.1788821)(5.283561,1.011515)
\psbezier[linewidth=0.04](2.9515827,-1.0)(1.7315826,-0.66)(0.30587173,0.37594724)(0.5515826,1.34)(0.7972935,2.3040528)(2.0875356,1.8207014)(2.6715827,1.24)(3.2556295,0.6592986)(3.2515826,0.44)(2.9915826,-0.96)
\psbezier[linewidth=0.04](4.2515826,0.12)(3.379995,0.76039475)(2.8458717,1.2996213)(2.914883,2.3998108)(2.983894,3.5)(5.981538,3.2553947)(6.0865602,2.0546052)(6.1915827,0.8538158)(5.0115824,0.32)(4.301179,0.16384251)
\psbezier[linewidth=0.04](2.4315827,-3.38)(1.2515826,-3.2)(0.0,-2.227836)(0.2515826,-1.06)(0.5031652,0.10783583)(2.2515826,-1.22)(2.9515827,-0.98)(3.6515825,-0.74)(3.9563887,-0.41792408)(4.3315825,-0.48)(4.7067766,-0.54207593)(5.1715827,-0.6)(5.3115826,-1.48)(5.4515824,-2.36)(3.3810325,-3.0947104)(2.4515827,-3.36)
\rput(1.0587701,1.05){\Large $C$} \rput(5.23877,2.19){\Large $C$}
\rput(3.5050201,-0.05){\Large $z$} \rput(4.614083,-1.69){\Large $L$}
\end{pspicture}
}
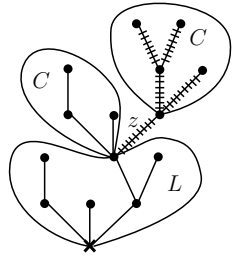
\\\vskip2mm\footnotesize{Figure 5. A right branch new type tree with a mutator.}
\end{center}

\begin{theorem}
The number of right branch new type trees with $n$ edges is
${2n-1\choose n}$. In particular, the proportion of new type
vertices is asymptotically $2\sqrt{{\pi\over n}}$.
\end{theorem}

\noindent\textbf{Proof.} The generating function for RBNT trees
follows from Figure 5:
\begin{eqnarray*}
T_M=L\cdot C\cdot zC={B\over C}\cdot C\cdot
zC=zBC={B-1\over2}=\sum_{n\ge1}{2n-1\choose n}z^n.
\end{eqnarray*}

The number of vertices involved are then counted by
\begin{eqnarray}\label{e:RBNT}
\left(zT_M\right)'={B-1\over2}+zB^3=\sum_{n\ge1}(n+1){2n-1\choose
n}z^n=2z+9z^2+40z^3+175z^4+\cdots
\end{eqnarray}
for (A097070). It also counts the number of parts equal to $1$ over
all weak compositions of $n+1$ into $n+1$ parts. Since the right
branch from the mutator contains all vertices of the new type, it
follows from Figure 5 that
\begin{eqnarray*}
V_N=L\cdot
C\cdot(z(zC))'=B\cdot(z^2C)'=B\cdot(2zC+z^2BC^2)=2+7z^2+26z^3+99z^4+\cdots.
\end{eqnarray*}
This is (A114121) in \cite{BSloa} except for the initial term. We
would like to find the proportion of new type vertices among all the
vertices. For $n\ge1$, the denominator is
$[z^n](zT_M)^\prime=(n+1){2n-1\choose n}={n+1\over2}{2n\choose n}$
and the numerator is
\begin{eqnarray*}
[z^n]V_N&=&[z^n](2zBC+z^2B^2C^2)=[z^n]\left(B-1+\left({B-1\over2}\right)^2\right)\\
&=&{2n\choose n}+{1\over4}\left(4^n-2{2n\choose
n}\right)=4^{n-1}+{1\over2}{2n\choose n}.
\end{eqnarray*}
Thus the proportion of new type vertices is
\begin{eqnarray*}
{4^{n-1}+{1\over2}{2n\choose n}\over{n+1\over2}{2n\choose
n}}\sim{\sqrt{\pi n}+{1\over2}\over {n+1\over2}}\sim2\sqrt{{\pi\over
n}}.
\end{eqnarray*}
\begin{center}
\begin{figure}
\end{figure}
\end{center}
\hfill{\rule{2mm}{2mm}}\bigskip

Another intermediate case of mutations in trees is obtained by
assuming that the vertices of the new type for trees go from the
mutator to the rightmost leaf above the mutator. For instance,
Figure 6 shows a tree with $5$ vertices of the new type. We call
trees involving such mutation the \textit{right path trees}.

\begin{figure}[h]
\begin{center}
\scalebox{0.6} {
\begin{pspicture}(0,-3.16)(4.1,3.14)
\psline[linewidth=0.04cm](3.4,1.76)(3.4,0.56)
\psdots[dotsize=0.24](3.4,1.78)
\psline[linewidth=0.04cm](3.28,1.6)(3.52,1.6)
\psline[linewidth=0.04cm](3.28,1.46)(3.52,1.46)
\psline[linewidth=0.04cm](3.28,1.3)(3.52,1.3)
\psline[linewidth=0.04cm](3.28,1.16)(3.52,1.16)
\psline[linewidth=0.04cm](3.28,1.04)(3.52,1.04)
\psline[linewidth=0.04cm](3.28,0.9)(3.52,0.9)
\psline[linewidth=0.04cm](3.28,0.76)(3.52,0.76)
\psline[linewidth=0.04cm](0.0,-0.52)(0.0,-1.72)
\psdots[dotsize=0.24](0.0,-0.52)
\psline[linewidth=0.08cm](1.46,-2.82)(1.74,-3.12)
\psline[linewidth=0.08cm](1.74,-2.82)(1.46,-3.12)
\psline[linewidth=0.04cm](1.6024584,-1.801012)(0.9972172,-0.483407)
\psdots[dotsize=0.24](1.0,-0.52)
\psline[linewidth=0.04cm](1.0,1.78)(1.0,0.58)
\psdots[dotsize=0.24](1.0,1.78)
\psline[linewidth=0.04cm](2.2,-0.62)(2.2,0.58)
\psline[linewidth=0.04cm](2.18,-0.56)(0.97721714,0.64)
\psdots[dotsize=0.24](1.0,0.58) \psdots[dotsize=0.24](2.2,0.58)
\psline[linewidth=0.04cm](2.2,-0.52)(3.4,0.58)
\psdots[dotsize=0.24](3.4,0.58)
\psline[linewidth=0.04cm](2.8832335,0.2999816)(3.055246,0.13261452)
\psline[linewidth=0.04cm](2.7856028,0.19964105)(2.9576151,0.032273963)
\psline[linewidth=0.04cm](2.6740248,0.08496613)(2.8460371,-0.082400955)
\psline[linewidth=0.04cm](2.576394,-0.015374423)(2.7484064,-0.18274151)
\psline[linewidth=0.04cm](2.4927104,-0.10138062)(2.664723,-0.2687477)
\psline[linewidth=0.04cm](2.3950796,-0.20172116)(2.567092,-0.36908826)
\psline[linewidth=0.04cm](2.2974489,-0.3020617)(2.4694612,-0.4694288)
\psline[linewidth=0.04cm](3.003817,0.37909964)(3.1758294,0.21173255)
\psline[linewidth=0.04cm](3.0999663,0.47888216)(3.2719786,0.31151506)
\psline[linewidth=0.04cm](3.4024584,1.7389879)(3.977208,3.0612826)
\psline[linewidth=0.04cm](3.4024584,1.7389879)(2.7972171,3.056593)
\psdots[dotsize=0.24](2.78,3.0) \psdots[dotsize=0.24](3.96,3.0)
\psline[linewidth=0.04cm](3.7398202,2.8549879)(3.9636145,2.7682922)
\psline[linewidth=0.04cm](3.6892476,2.724441)(3.913042,2.6377456)
\psline[linewidth=0.04cm](3.6314507,2.575245)(3.8552449,2.4885492)
\psline[linewidth=0.04cm](3.5808783,2.4446983)(3.8046725,2.3580027)
\psline[linewidth=0.04cm](3.5375304,2.332801)(3.7613246,2.2461054)
\psline[linewidth=0.04cm](3.486958,2.2022545)(3.7107522,2.1155589)
\psline[linewidth=0.04cm](3.4363856,2.0717077)(3.6601799,1.9850123)
\psline[linewidth=0.04cm](1.6,-2.92)(1.6,-1.72)
\psline[linewidth=0.04cm](1.6,-2.92)(0.0,-1.72)
\psdots[dotsize=0.24](0.0,-1.72) \psdots[dotsize=0.24](1.6,-1.76)
\psline[linewidth=0.04cm](1.6,-1.82)(2.2,-0.52)
\psdots[dotsize=0.24](2.2,-0.52)
\psline[linewidth=0.04cm](1.9398202,-0.70501226)(2.1636145,-0.7917078)
\psline[linewidth=0.04cm](1.8892478,-0.8355589)(2.1130419,-0.9222545)
\psline[linewidth=0.04cm](1.8314507,-0.9847551)(2.055245,-1.0714507)
\psline[linewidth=0.04cm](1.7808782,-1.1153017)(2.0046725,-1.2019973)
\psline[linewidth=0.04cm](1.7375305,-1.2271988)(1.9613247,-1.3138945)
\psline[linewidth=0.04cm](1.686958,-1.3577455)(1.9107523,-1.4444411)
\psline[linewidth=0.04cm](1.6363856,-1.4882922)(1.8601799,-1.5749878)
\pscircle[linewidth=0.03,dimen=outer](1.6,-1.76){0.26}
\psline[linewidth=0.04cm](3.2,-0.52)(3.2,-1.72)
\psdots[dotsize=0.24](3.2,-0.52)
\psline[linewidth=0.04cm](1.6,-2.92)(3.2,-1.72)
\psdots[dotsize=0.24](3.2,-1.72)
\end{pspicture}
}\\\vskip2mm\footnotesize{Figure 6. A right path tree with $5$ new
type vertices.}
\end{center}
\end{figure}

\begin{theorem}
The number of right path trees with $n$ edges is ${2n\choose n}$. In
particular, the proportion of new type vertices is asymptotically
${1\over2n}$.
\end{theorem}

\noindent\textbf{Proof.} It can be easily seen that the generating
function for right path trees is
\begin{eqnarray*}
T_M=L\cdot(1+zC+z^2C^2+\cdots)={B\over C}\cdot{1\over1-zC}={B\over
C}\cdot C=B.
\end{eqnarray*}
Hence the number of vertices in right path trees is counted by
$(zB)'=B+2zB^3=\sum_{n\ge0}(n+1){2n\choose
n}z^n=1+4z+18z^2+80z^3+\cdots$. Next we want to compute the number
of new type vertices. If the right path from the mutator has length
$k$, we have $k+1$ new type vertices and the generating function for
vertices $(k+1)z^kC^k$ since a subtree can be attached to the left
of all vertices except the last. Including the location of the
mutator and summing over $k$ gives the generating function
\begin{eqnarray*}
V_N=L\cdot(1+2zC+3z^2C^2+\cdots)={B\over
C}\cdot{1\over(1-zC)^2}={B\over C}\cdot
C^2=BC=\sum_{n\ge0}{2n+1\choose n}z^n.
\end{eqnarray*}
Thus the proportion of new type vertices is
\begin{eqnarray*}
{{1\over2}{2n\choose n}\over (n+1){2n\choose
n}}={1\over2(n+1)}\sim{1\over2n}.
\end{eqnarray*}
\hfill{\rule{2mm}{2mm}}\bigskip

A common situation is that the mutator is not known until the child
appears in the new state. So we want to look at this rightmost path
from the mutator to the new type leaf but now requiring the mutator
to have at least one descendant. Such mutation will be called the
\textit{right path$^\ast$} mutation.

The number of right path$^\ast$ trees has the generating function
\begin{eqnarray*}
L\cdot(zC+z^2C^2+\cdots)={B\over
C}\cdot{zC\over1-zC}=zB\cdot{1\over1-zC}=zBC=\sum_{n\ge1}{2n-1\choose
n}z^n,
\end{eqnarray*}
which is the same as $T_M$ of right branch new type trees. So the
generating function for vertices is given by \eref{e:RBNT}.

For the number of vertices of the new type we have
\begin{align*}
L\cdot(2zC+3(zC)^2+\cdots)&={B\over
C}\cdot\left({1\over(1-zC)^2}-1\right)={B\over
C}\cdot(C^2-1)=BC-{B\over C}\\&=\sum_{n\ge1}{3n+1\over
2n+2}{2n\choose n}z^n =2z+7z^2+25z^3+91z^4+\cdots.
\end{align*}
This is (A097613) in \cite{BSloa} except for the initial term. The
proportion of new type vertices out of all vertices is
\begin{eqnarray*}
{{3n+1\over 2n+2}{2n\choose n}\over {n+1\over2}{2n\choose
n}}={3n+1\over(n+1)^2}\sim{3\over n}.
\end{eqnarray*}
Here is an illustration for the $10$ right path$^\ast$ trees on $3$
edges. There are $40$ vertices of which $25$ are of the new type.
\begin{figure}[h]
\begin{center}
\scalebox{0.6} {
\begin{pspicture}(0,-1.99)(16.80219,2.01)
\psline[linewidth=0.04cm](0.26,1.85)(0.26,-1.75)
\psline[linewidth=0.04cm](7.6824584,-1.791012)(7.077217,-0.473407)
\psdots[dotsize=0.24](0.26,1.87) \psdots[dotsize=0.24](0.26,0.67)
\psline[linewidth=0.08cm](0.12,-1.57)(0.4,-1.87)
\psline[linewidth=0.08cm](0.4,-1.57)(0.12,-1.87)
\psline[linewidth=0.04cm](3.68,-0.53)(3.68,-1.73)
\psdots[dotsize=0.24](3.68,-0.53)
\psline[linewidth=0.08cm](3.52,-1.55)(3.8,-1.85)
\psline[linewidth=0.08cm](3.8,-1.55)(3.52,-1.85)
\psdots[dotsize=0.24](7.06,-0.53)
\psline[linewidth=0.08cm](7.52,-1.57)(7.8,-1.87)
\psline[linewidth=0.08cm](7.8,-1.57)(7.52,-1.87)
\pscircle[linewidth=0.03,dimen=outer](0.26,-1.73){0.26}
\pscircle[linewidth=0.03,dimen=outer](3.66,-1.69){0.26}
\psdots[dotsize=0.24](0.26,-0.53)
\pscircle[linewidth=0.03,dimen=outer](7.66,-1.71){0.26}
\psline[linewidth=0.04cm](7.04,0.63)(7.04,-0.57)
\psdots[dotsize=0.24](7.06,0.67)
\psline[linewidth=0.04cm](3.7024584,-0.59101206)(3.097217,0.726593)
\psdots[dotsize=0.24](3.08,0.67)
\psline[linewidth=0.04cm](11.722459,-1.791012)(11.117217,-0.473407)
\psdots[dotsize=0.24](11.1,-0.53)
\psline[linewidth=0.08cm](11.56,-1.57)(11.84,-1.87)
\psline[linewidth=0.08cm](11.84,-1.57)(11.56,-1.87)
\pscircle[linewidth=0.03,dimen=outer](11.7,-1.73){0.26}
\psline[linewidth=0.04cm](12.26,0.61)(12.26,-0.59)
\psdots[dotsize=0.24](12.28,0.65)
\psline[linewidth=0.04cm](1.26,1.85)(1.26,-1.75)
\psdots[dotsize=0.24](1.26,1.87) \psdots[dotsize=0.24](1.26,0.67)
\psline[linewidth=0.08cm](1.12,-1.57)(1.4,-1.87)
\psline[linewidth=0.08cm](1.4,-1.57)(1.12,-1.87)
\pscircle[linewidth=0.03,dimen=outer](1.26,-0.53){0.26}
\psdots[dotsize=0.24](1.26,-0.53)
\psline[linewidth=0.04cm](5.68,-0.55)(5.68,-1.75)
\psdots[dotsize=0.24](5.68,-0.55)
\psline[linewidth=0.08cm](5.52,-1.57)(5.8,-1.87)
\psline[linewidth=0.08cm](5.8,-1.57)(5.52,-1.87)
\pscircle[linewidth=0.03,dimen=outer](5.66,-0.55){0.26}
\psline[linewidth=0.04cm](5.7024584,-0.61101204)(5.097217,0.706593)
\psdots[dotsize=0.24](5.08,0.65)
\psline[linewidth=0.04cm](9.682459,-1.7710121)(10.257208,-0.44871742)
\psline[linewidth=0.04cm](9.682459,-1.7710121)(9.077217,-0.453407)
\psdots[dotsize=0.24](9.06,-0.51) \psdots[dotsize=0.24](10.24,-0.51)
\psline[linewidth=0.08cm](9.52,-1.55)(9.8,-1.85)
\psline[linewidth=0.08cm](9.8,-1.55)(9.52,-1.85)
\pscircle[linewidth=0.03,dimen=outer](9.04,-0.51){0.26}
\psline[linewidth=0.04cm](9.04,0.65)(9.04,-0.55)
\psdots[dotsize=0.24](9.06,0.69)
\psline[linewidth=0.04cm](13.682459,-1.791012)(14.257208,-0.4687174)
\psline[linewidth=0.04cm](13.682459,-1.791012)(13.077217,-0.473407)
\psdots[dotsize=0.24](13.06,-0.53)
\psdots[dotsize=0.24](14.24,-0.53)
\psline[linewidth=0.08cm](13.52,-1.57)(13.8,-1.87)
\psline[linewidth=0.08cm](13.8,-1.57)(13.52,-1.87)
\pscircle[linewidth=0.03,dimen=outer](14.22,-0.55){0.26}
\psline[linewidth=0.04cm](14.22,0.61)(14.22,-0.59)
\psdots[dotsize=0.24](14.24,0.65)
\psline[linewidth=0.04cm](15.86,-1.711012)(15.86,-0.51)
\psdots[dotsize=0.24](15.86,-0.51)
\psline[linewidth=0.08cm](15.72,-1.57)(16.0,-1.87)
\psline[linewidth=0.08cm](16.0,-1.57)(15.72,-1.87)
\pscircle[linewidth=0.03,dimen=outer](15.86,-1.73){0.26}
\psline[linewidth=0.04cm](15.882459,-1.7710121)(15.08,-0.59)
\psdots[dotsize=0.24](15.06,-0.53)
\psline[linewidth=0.04cm](2.26,1.85)(2.26,-1.75)
\psdots[dotsize=0.24](2.26,1.87) \psdots[dotsize=0.24](2.26,0.67)
\psline[linewidth=0.08cm](2.12,-1.57)(2.4,-1.87)
\psline[linewidth=0.08cm](2.4,-1.57)(2.12,-1.87)
\pscircle[linewidth=0.03,dimen=outer](2.26,0.67){0.26}
\psdots[dotsize=0.24](2.26,-0.53)
\psline[linewidth=0.04cm](0.14,-0.73)(0.38,-0.73)
\psline[linewidth=0.04cm](0.14,-0.87)(0.38,-0.87)
\psline[linewidth=0.04cm](0.14,-1.03)(0.38,-1.03)
\psline[linewidth=0.04cm](0.14,-1.17)(0.38,-1.17)
\psline[linewidth=0.04cm](0.14,-1.29)(0.38,-1.29)
\psline[linewidth=0.04cm](0.14,-1.43)(0.38,-1.43)
\psline[linewidth=0.04cm](0.14,-1.57)(0.38,-1.57)
\psline[linewidth=0.04cm](0.14,0.49)(0.38,0.49)
\psline[linewidth=0.04cm](0.14,0.35)(0.38,0.35)
\psline[linewidth=0.04cm](0.14,0.19)(0.38,0.19)
\psline[linewidth=0.04cm](0.14,0.05)(0.38,0.05)
\psline[linewidth=0.04cm](0.14,-0.07)(0.38,-0.07)
\psline[linewidth=0.04cm](0.14,-0.21)(0.38,-0.21)
\psline[linewidth=0.04cm](0.14,-0.35)(0.38,-0.35)
\psline[linewidth=0.04cm](0.14,1.69)(0.38,1.69)
\psline[linewidth=0.04cm](0.14,1.55)(0.38,1.55)
\psline[linewidth=0.04cm](0.14,1.39)(0.38,1.39)
\psline[linewidth=0.04cm](0.14,1.25)(0.38,1.25)
\psline[linewidth=0.04cm](0.14,1.13)(0.38,1.13)
\psline[linewidth=0.04cm](0.14,0.99)(0.38,0.99)
\psline[linewidth=0.04cm](0.14,0.85)(0.38,0.85)
\psline[linewidth=0.04cm](1.14,0.41)(1.38,0.41)
\psline[linewidth=0.04cm](1.14,0.27)(1.38,0.27)
\psline[linewidth=0.04cm](1.14,0.11)(1.38,0.11)
\psline[linewidth=0.04cm](1.14,-0.03)(1.38,-0.03)
\psline[linewidth=0.04cm](1.14,-0.15)(1.38,-0.15)
\psline[linewidth=0.04cm](1.14,-0.29)(1.38,-0.29)
\psline[linewidth=0.04cm](1.14,-0.43)(1.38,-0.43)
\psline[linewidth=0.04cm](1.14,1.69)(1.38,1.69)
\psline[linewidth=0.04cm](1.14,1.55)(1.38,1.55)
\psline[linewidth=0.04cm](1.14,1.39)(1.38,1.39)
\psline[linewidth=0.04cm](1.14,1.25)(1.38,1.25)
\psline[linewidth=0.04cm](1.14,1.13)(1.38,1.13)
\psline[linewidth=0.04cm](1.14,0.99)(1.38,0.99)
\psline[linewidth=0.04cm](1.14,0.85)(1.38,0.85)
\psline[linewidth=0.04cm](2.14,1.71)(2.38,1.71)
\psline[linewidth=0.04cm](2.14,1.57)(2.38,1.57)
\psline[linewidth=0.04cm](2.14,1.41)(2.38,1.41)
\psline[linewidth=0.04cm](2.14,1.27)(2.38,1.27)
\psline[linewidth=0.04cm](2.14,1.15)(2.38,1.15)
\psline[linewidth=0.04cm](2.14,1.01)(2.38,1.01)
\psline[linewidth=0.04cm](2.14,0.87)(2.38,0.87)
\psline[linewidth=0.04cm](3.56,-0.73)(3.8,-0.73)
\psline[linewidth=0.04cm](3.56,-0.87)(3.8,-0.87)
\psline[linewidth=0.04cm](3.56,-1.03)(3.8,-1.03)
\psline[linewidth=0.04cm](3.56,-1.17)(3.8,-1.17)
\psline[linewidth=0.04cm](3.56,-1.29)(3.8,-1.29)
\psline[linewidth=0.04cm](3.56,-1.43)(3.8,-1.43)
\psline[linewidth=0.04cm](3.56,-1.57)(3.8,-1.57)
\psline[linewidth=0.04cm](5.6624584,-0.61101204)(6.237208,0.7112826)
\psdots[dotsize=0.24](6.22,0.65)
\psline[linewidth=0.04cm](5.99982,0.5049878)(6.223614,0.41829216)
\psline[linewidth=0.04cm](5.949248,0.37444112)(6.173042,0.2877455)
\psline[linewidth=0.04cm](5.891451,0.22524492)(6.115245,0.13854934)
\psline[linewidth=0.04cm](5.840878,0.09469827)(6.0646725,0.008002682)
\psline[linewidth=0.04cm](5.7975307,-0.01719886)(6.0213246,-0.10389445)
\psline[linewidth=0.04cm](5.746958,-0.14774552)(5.9707522,-0.2344411)
\psline[linewidth=0.04cm](5.6963854,-0.27829218)(5.92018,-0.36498776)
\psline[linewidth=0.04cm](7.6624584,-1.791012)(8.237207,-0.4687174)
\psdots[dotsize=0.24](8.22,-0.53)
\psline[linewidth=0.04cm](7.99982,-0.67501223)(8.223615,-0.76170784)
\psline[linewidth=0.04cm](7.949248,-0.8055589)(8.173042,-0.8922545)
\psline[linewidth=0.04cm](7.891451,-0.95475507)(8.115245,-1.0414506)
\psline[linewidth=0.04cm](7.840878,-1.0853018)(8.064672,-1.1719973)
\psline[linewidth=0.04cm](7.7975307,-1.1971989)(8.021325,-1.2838944)
\psline[linewidth=0.04cm](7.746958,-1.3277456)(7.9707522,-1.4144411)
\psline[linewidth=0.04cm](7.6963854,-1.4582921)(7.92018,-1.5449878)
\psline[linewidth=0.04cm](15.899038,-1.686968)(16.679945,-0.4749503)
\psdots[dotsize=0.24,dotangle=-9.301077](16.653057,-0.53264606)
\psline[linewidth=0.04cm](16.412336,-0.64016575)(16.619175,-0.7618917)
\psline[linewidth=0.04cm](16.34133,-0.7608224)(16.548168,-0.8825484)
\psline[linewidth=0.04cm](16.260178,-0.8987158)(16.467018,-1.0204418)
\psline[linewidth=0.04cm](16.18917,-1.0193725)(16.396011,-1.1410984)
\psline[linewidth=0.04cm](16.12831,-1.1227925)(16.335148,-1.2445184)
\psline[linewidth=0.04cm](16.057302,-1.2434492)(16.264141,-1.3651751)
\psline[linewidth=0.04cm](15.986295,-1.3641058)(16.193134,-1.4858317)
\psline[linewidth=0.04cm](3.6824584,-0.59101206)(4.257208,0.7312826)
\psdots[dotsize=0.24](4.24,0.67)
\psline[linewidth=0.04cm](4.01982,0.52498776)(4.243614,0.43829218)
\psline[linewidth=0.04cm](3.9692478,0.3944411)(4.193042,0.30774552)
\psline[linewidth=0.04cm](3.9114506,0.24524494)(4.135245,0.15854934)
\psline[linewidth=0.04cm](3.8608782,0.114698276)(4.0846725,0.028002681)
\psline[linewidth=0.04cm](3.8175304,0.0028011396)(4.0413246,-0.083894454)
\psline[linewidth=0.04cm](3.766958,-0.12774551)(3.9907522,-0.2144411)
\psline[linewidth=0.04cm](3.7163856,-0.25829217)(3.9401798,-0.34498775)
\psline[linewidth=0.04cm](8.92,0.51)(9.16,0.51)
\psline[linewidth=0.04cm](8.92,0.37)(9.16,0.37)
\psline[linewidth=0.04cm](8.92,0.21)(9.16,0.21)
\psline[linewidth=0.04cm](8.92,0.07)(9.16,0.07)
\psline[linewidth=0.04cm](8.92,-0.05)(9.16,-0.05)
\psline[linewidth=0.04cm](8.92,-0.19)(9.16,-0.19)
\psline[linewidth=0.04cm](8.92,-0.33)(9.16,-0.33)
\psline[linewidth=0.04cm](12.14,0.49)(12.38,0.49)
\psline[linewidth=0.04cm](12.14,0.35)(12.38,0.35)
\psline[linewidth=0.04cm](12.14,0.19)(12.38,0.19)
\psline[linewidth=0.04cm](12.14,0.05)(12.38,0.05)
\psline[linewidth=0.04cm](12.14,-0.07)(12.38,-0.07)
\psline[linewidth=0.04cm](12.14,-0.21)(12.38,-0.21)
\psline[linewidth=0.04cm](12.14,-0.35)(12.38,-0.35)
\psline[linewidth=0.04cm](14.1,0.47)(14.34,0.47)
\psline[linewidth=0.04cm](14.1,0.33)(14.34,0.33)
\psline[linewidth=0.04cm](14.1,0.17)(14.34,0.17)
\psline[linewidth=0.04cm](14.1,0.03)(14.34,0.03)
\psline[linewidth=0.04cm](14.1,-0.09)(14.34,-0.09)
\psline[linewidth=0.04cm](14.1,-0.23)(14.34,-0.23)
\psline[linewidth=0.04cm](14.1,-0.37)(14.34,-0.37)
\psline[linewidth=0.04cm](11.722459,-1.791012)(12.297208,-0.4687174)
\psdots[dotsize=0.24](12.28,-0.53)
\psline[linewidth=0.04cm](12.05982,-0.67501223)(12.283614,-0.76170784)
\psline[linewidth=0.04cm](12.009248,-0.8055589)(12.233042,-0.8922545)
\psline[linewidth=0.04cm](11.95145,-0.95475507)(12.175245,-1.0414506)
\psline[linewidth=0.04cm](11.900878,-1.0853018)(12.124673,-1.1719973)
\psline[linewidth=0.04cm](11.857531,-1.1971989)(12.081325,-1.2838944)
\psline[linewidth=0.04cm](11.806958,-1.3277456)(12.030752,-1.4144411)
\psline[linewidth=0.04cm](11.756386,-1.4582921)(11.98018,-1.5449878)
\end{pspicture}
}\\\vskip2mm\footnotesize{Figure 7. Right path trees with a mutator
having at least one child.}
\end{center}
\end{figure}

The following table summarizes the main results in this section.
\begin{center}
\begin{tabular}{l||*{1}{l}|*{1}{l}}\hline
Mutation&Number of trees&Asymptotic ratio of\\
&with the mutation&vertices of the new
type\\\specialrule{0.5mm}{0.5mm}{0.5mm} Short lived&${2n-2\choose n-1}$&$2/((n+1){2n-2\choose n-1})$\\
Toggle&${2n+1\choose n}$&${1\over2}\sqrt{{\pi\over n}}$\\
Embedded new type&$\sum_{k=0}^n{1\over k+1}{2n\choose n-k}{2k\choose k}$&${3\over 5}$\\
Right branch new type&${2n-1\choose n}$&$2\sqrt{{\pi\over n}}$\\
Right path&${2n\choose n}$&${1\over 2n}$\\
Right path$^\ast$&${2n-1\choose n}$&${3\over n}$\\ \hline
\end{tabular}
\end{center}

\section{Another class of ordered trees}

What if a mutation occurs in another class of ordered trees instead
of the usual ordered trees? We end this paper with the complete
binary trees with a simple mutation in which a mutator changes all
of its descendants to the new type. A complete binary tree has every
vertex with updegree 0 or 2. It is known that there are $C_n$
complete binary trees with $n$ internal vertices. To have the
coefficient of $z^n$ count edges instead of internal vertices, the
appropriate generating function is
$\sum_{n\ge0}C_nz^{2n}=C(z^2)={1-\sqrt{1-4z^2}\over2z^2}$. This
change brings about various other small changes. By setting
$\tilde{C}(z)=\tilde{C}=C(z^2)$ and $\tilde{B}(z)=\tilde{B}=B(z^2)$,
we have $\tilde{C}'=2z\tilde{B}\tilde{C}^2$,
$\tilde{B}'=4z\tilde{B}^3$ and $\tilde{B}=1+2z^2\tilde{B}\tilde{C}$.
The vertex and leaf generating functions $V$ and $L$ are
\begin{eqnarray*}
V=(z\tilde{C})'=\tilde{C}\cdot(1+2z^2\tilde{B}\tilde{C})
=\tilde{C}\tilde{B}\quad\hbox{and}\quad
L={\tilde{C}\tilde{B}\over\tilde{C}}=\tilde{B}.
\end{eqnarray*}

Let $T_M$ be the generating function for complete binary trees with
a mutator, where all the descendants of the mutator are of the new
type. Then $T_M=V=\tilde{C}\tilde{B}=1+3z^2+10z^4+35z^6+\cdots$. The
number of vertices of these mutator enhanced trees has the
generating function
\begin{align*}
\tilde{V}&=(z\tilde{T})'=\tilde{T}+z\tilde{T}'=\tilde{B}\tilde{C}+z(\tilde{B}\tilde{C})'
=\tilde{B}\tilde{C}+z(\tilde{B}\cdot2z\tilde{B}\tilde{C}^2+4z\tilde{B}^3\tilde{C})\\
&=\tilde{B}\tilde{C}\cdot(1+2z^2\tilde{B}\tilde{C}+4z^2\tilde{B}^2)
=\tilde{B}\tilde{C}\cdot\left(\tilde{B}+{4z^2\over1-4z^2}\right)\\&=\sum_{n\ge0}(2n+1)^2C_nz^{2n}
=1+9z^2+50z^4+245z^6+\cdots.
\end{align*}
This is a new entry in the OEIS \cite{BSloa}.

To figure out the number of new type vertices, we apply the uplift
principle again and get the generating function
\begin{eqnarray*}
L\tilde{V}=\tilde{B}\cdot\tilde{B}\tilde{C}=\tilde{B}^2\tilde{C}=
\sum_{n\ge0}\left(2^{2n-1}-{1\over2}{2n\choose n}\right)z^{2n}
=1+5z^2+22z^4+93z^6+\cdots.
\end{eqnarray*}
The proportion of new type vertices is
\begin{eqnarray*}
{[z^n]\tilde{B}^2\tilde{C}\over[z^n]\tilde{V}}
={2^{2n-1}-{1\over2}{2n\choose n}\over(2n+1)^2{1\over n+1}{2n\choose
n}}\sim {{1\over2}{2n\choose n}\sqrt{\pi n}-{1\over2}{2n\choose
n}\over {(2n+1)^2\over n+1}{2n\choose n}}= {{1\over2}\sqrt{\pi
n}-{1\over2}\over {(2n+1)^2\over n+1}}\sim{1\over8}\sqrt{{\pi\over
n}}.
\end{eqnarray*}

Complete binary trees and ordered trees are both examples of trees
satisfying a uniform updegree requirement. The methods used in this
paper generalize to all such trees but that will be the subject of a
separate manuscript.

\end{document}